\documentclass[12pt]{article}
\usepackage{amsmath,amssymb}
\usepackage[dvipdfmx]{graphicx}

\begin{document}

\title{Divergent series in Gauss's diary and their extensions}

\author{Kiyoshi \textsc{Sogo}
\thanks{EMail: sogo@icfd.co.jp}
}

\date{}

\maketitle

\begin{center}
Institute of Computational Fluid Dynamics, 1-16-5, Haramachi, Meguroku, 
Tokyo, 152-0011, Japan
\end{center}

\abstract{
Two divergent series in Entry 7 of Gauss's diary are extended systematically by introducing additional parameters. 
Rogers-Fine identities, Ramanujan's continued fractions and Heine's transformation relations of 
basic hypergeometric series are applied to find equivalent alternative series, which enable us to compute 
the sums of divergent series in question. 
}
\tableofcontents

\section{Introduction}
\setcounter{equation}{0}

In 1796 Gauss wrote, in Entry 7 of his mathematical diary \cite{GaussD, Dunnington}, two formulas such that
\begin{align}
& 1-2+8-64+\cdots
\nonumber \\
&\qquad=\frac{1}{1+}\ \frac{2}{1+}\ \frac{2}{1+}\ \frac{8}{1+}\ \frac{12}{1+}\ \frac{32}{1+}\ 
\frac{56}{1+}\ \frac{128}{1+}\ \cdots,
\label{Gauss01}
\\
& 1-1+1\cdot 3-1\cdot 3\cdot 7+1\cdot 3\cdot 7\cdot 15-\cdots
\nonumber \\
&\qquad=\frac{1}{1+}\ \frac{1}{1+}\ \frac{2}{1+}\ \frac{6}{1+}\ \frac{12}{1+}\ \frac{28}{1+}\ \cdots,
\label{Gauss02}
\end{align}
where the right hand sides are continued fractions, which we use the notation 
\begin{align}
\frac{d_0}{1+}\ \frac{d_1}{1+}\ \frac{d_2}{1+}\ \cdots=\frac{d_0}{1+\frac{d_1}{1+\frac{d_2}{1+\cdots}}}.
\end{align}
They are the $q=2$ case of 
\begin{align}
&\sum_{n=0}^\infty (-1)^n q^{n(n+1)/2}=
\frac{1}{1+}\ \frac{q}{1+}\ \frac{q^2-q}{1+}\ \frac{q^3}{1+}\ \frac{q^4-q^2}{1+}\ \cdots,
\label{Gauss1}
\\
&\sum_{n=0}^\infty (q;q)_n=
\frac{1}{1+}\ \frac{q-1}{1+}\ \frac{q^2-q}{1+}\ \frac{q^3-q}{1+}\ \frac{q^4-q^2}{1+}\ \cdots,
\label{Gauss2}
\end{align}
respectively, where $(\alpha;q)_n$ is $q$-Pochhammer symbol defined by
\begin{align}
(\alpha;q)_0=1,\ &(\alpha; q)_n=(\alpha;q)_{n-1}\cdot(1-\alpha q^{n-1}),\quad (n\geq 1)
\end{align}
which gives $(q;q)_0=1,\ (q;q)_n=(1-q)\cdots(1-q^n),\ (n\geq 1)$. 
For the simplicity, the parameter $q$ is assumed for a while to be real and non-negative number. 

Since the left hand side of \eqref{Gauss1} and \eqref{Gauss2} are divergent when $q>1$, we can suppose that here Gauss proposed 
problems to find {\it sums of divergent series} similarly as Euler did before \cite{Bourbaki, HardyD}. 
The purpose of the present article is to give answers to \eqref{Gauss01} and \eqref{Gauss02}. 
In addition we wish to propose their extensions such that 
\begin{align}
S_\rho^{(1)}(q)&=\sum_{n=0}^\infty (-1)^n q^{n[(\rho-2)n-(\rho-4)]/2},\quad (\rho=3,4,5,\cdots),
\label{G3}
\\
&=1-q+q^\rho-q^{3\rho-3}+\cdots,
\nonumber
\\
S_\kappa^{(2)}(q)&=\sum_{n=0}^\infty q^n (q;q^\kappa)_n,\quad(\kappa=1,2,3,\cdots),
\label{G4}
\\
&=1+q(1-q)+q^2(1-q)(1-q^{1+\kappa})+\cdots,
\nonumber
\end{align}
and give values of them for $q>1$, as well as for $q<1$. 
It should be noted that the power $n[(\rho-2)n-(\rho-4)]/2$ in \eqref{G3} is called $\rho$ {\it polygonal number} series, 
that is, $\rho=3$ triangle, $\rho=4$ square, $\rho=5$ pentagon etc., 
and the first case of Gauss \eqref{Gauss1} corresponds to $\rho=3$ trigonal number series. 

We should also note that \eqref{G4} is $x=q$ case of the power series
\begin{align}
S_\kappa^{(2)}(x)&=\sum_{n=0}^\infty (q; q^\kappa)_n x^n
\label{P_CF}
\\
&=\frac{1}{1-}\ \frac{(1-q)x}{1-}\ \frac{q(1-q^\kappa)x}{1-}\ \frac{q^\kappa(1-q^{1+\kappa})x}{1-}\ 
\frac{q^{1+\kappa}(1-q^{2\kappa})x}{1-}\ \cdots,
\nonumber
\end{align}
and the second case of Gauss \eqref{Gauss2} is given by setting $\kappa=1$ and $x=1$, {\it i.e.}, $S_1^{(2)}(1)$. 
The continued fraction expression in the right hand side will be derived later. 


\section{Answers to Gauss problems}
\setcounter{equation}{0}


\subsection{First group of problems -- polygonal number series}

By using Rogers-Fine identity \cite{RogersFine1, RogersFine2, McLaurin}
\begin{align}
\sum_{n=0}^\infty x^n q^{n(n+1)/2} = \sum_{n=0}^\infty \frac{(xq; q^2)_n}{(xq^2; q^2)_n}\cdot (xq)^n,
\label{RogersFine}
\end{align}
at $x=-1$, we obtain 
\begin{align}
\sum_{n=0}^\infty (-1)^n q^{n(n+1)/2}=1+\sum_{n=1}^\infty (-1)^n q^n\cdot
\frac{(1+q)\cdots(1+q^{2n-1})}{(1+q^2)\cdots(1+q^{2n})}.
\label{Gauss1_equality_q}
\end{align}
Obviously both sides are convergent for $q<1$. But for $q>1$, 
although the left hand side is divergent, the right hand side becomes by substituting $q=p^{-1}$ 
\begin{align}
&1+\sum_{n=1}^\infty (-1)^n q^n\cdot \frac{(1+q)\cdots(1+q^{2n-1})}{(1+q^2)\cdots(1+q^{2n})}
\nonumber \\
&\qquad\qquad =
1+\sum_{n=1}^\infty (-1)^n \frac{(1+p)\cdots(1+p^{2n-1})}{(1+p^2)\cdots(1+p^{2n})}.
\label{Gauss1_equality_p}
\end{align}
The right hand side of \eqref{Gauss1_equality_p} is {\it almost convergent} when $p<1$, {\it i.e.}, $q>1$, 
because it is similar to the alternating series of $Q_\infty$
\begin{align}
\sum_{n=0}^\infty (-1)^n Q_\infty,\qquad Q_\infty=\prod_{n=1}^\infty \left(\frac{1+p^{2n-1}}{1+p^{2n}}\right),
\label{Qinfinity}
\end{align}
whose sum is given by $\frac{1}{2}Q_\infty$ in the sense of Ces\`{a}ro sum \cite{HardyD}.
Therefore by using \eqref{Gauss1_equality_p} we can {\it define and replace} the left hand side for $q>1$ 
by the right hand side for $p=q^{-1}<1$. 

Result of numerical experiment of partial sums for $p=0.5\ (q=2)$ is given by
\begin{align}
S_N&= 1+\sum_{n=1}^N (-1)^n \frac{(1+p)\cdots(1+p^{2n-1})}{(1+p^2)\cdots(1+p^{2n})}
\\
\begin{split}
&= +1.0759457568\cdots\quad(N=100),\\
&= -0.2208954963\cdots\quad(N=101),
\end{split}
\end{align}
repeating these two values alternately as $N$ increases.  
The above values are convergent already for $N=50$ up to the given ten digits. 
Since $Q_\infty\simeq S_{N}-S_{N+1}$ for large even $N$, we can check the equality 
\begin{align}
Q_\infty=1.296841253\cdots=1.0759457568\cdots+0.2208954963\cdots
\end{align}
by computing $Q_\infty$ independently by \eqref{Qinfinity}. 
In other words, we can say for large even $N$ 
\begin{align}
S_N\simeq Lim+\frac{1}{2} Q_\infty,\quad
S_{N+1}\simeq Lim-\frac{1}{2} Q_\infty,
\end{align}
and the value of $Lim$ is given by the average for $p=0.5,\ N=100$,
\begin{align}
1-2+2^3-2^6+\cdots&=Lim
\nonumber \\
&\simeq\frac{1}{2}\left(S_N+S_{N+1}\right)=0.4275251302\cdots,
\label{GaussP1}
\end{align}
which is our answer to Gauss problem 1, {\it i.e.}, \eqref{Gauss01}. \\

Now let us consider some extensions of this first problem in Entry 7 of Gauss's diary. 
If we change $q\rightarrow q^a,\ x\rightarrow xq^{-b}$ in \eqref{RogersFine}, we have
\begin{align}
\sum_{n=0}^\infty x^n q^{an(n+1)/2-nb}
=\sum_{n=0}^\infty x^n q^{(a-b)n}\ \frac{(xq^{a-b}; q^{2a})_n}{(xq^{2a-b}; q^{2a})_n}.
\end{align}
By setting $a=\rho-2,\ b=\rho-3$ and $x=-1$, the above equality becomes
\begin{align}
\begin{split}
S_\rho^{(1)}(q)&=\sum_{n=0}^\infty (-1)^n q^{n[(\rho-2)n-(\rho-4)]/2}
\\
&=\sum_{n=0}^\infty (-1)^n q^n\ \frac{(-q; q^{2\rho-4})_n}{(-q^{\rho-1}; q^{2\rho-4})_n}
\\
&=\sum_{n=0}^\infty (-1)^n p^{(\rho-3)n} \frac{(-p;p^{2\rho-4})_n}{(-p^{\rho-1}; p^{2\rho-4})_n},
\end{split}
\end{align}
where $p=q^{-1}$ as before. In this way, we arrive at the following theorem. \\

\noindent
{\bf Theorem 1}\quad For $\rho=3,4,5,\cdots$, we have ($p=q^{-1}$)
\begin{align}
&S_\rho^{(1)}(q)=\sum_{n=0}^\infty (-1)^n q^{n[(\rho-2) n-(\rho-4)]/2}
\nonumber \\
&\quad=1+\sum_{n=1}^\infty (-1)^n q^n\frac{(1+q)\cdots(1+q^{1+(2\rho-4)(n-1)})}
{(1+q^{\rho-1})\cdots(1+q^{\rho-1+(2\rho-4)(n-1)})}
\label{GF_1}
\\
&\quad=1+\sum_{n=1}^\infty (-1)^n p^{(\rho-3)n} \frac{(1+p)\cdots(1+p^{1+(2\rho-4)(n-1)})}
{(1+p^{\rho-1})\cdots(1+p^{\rho-1+(2\rho-4)(n-1)})},
\nonumber
\end{align}
which gives \eqref{Gauss1_equality_q},\ \eqref{Gauss1_equality_p} when $\rho=3$, as it should be.\\

Let us consider $q$ as a complex number. When $\rho>3$, the first line and the second line in \eqref{GF_1} are convergent 
for $|q|<1$ and coincide with each other. 
Now the second line is convergent even for $|q|>1$, which is equivalent to the third line for $|p|<1$, because the third line is 
just a rewriting of the second line by setting $q=p^{-1}$. And since the third line is convergent for $|p|<1$ {\it i.e.}, $|q|>1$, 
we can say the third line is {\it analytical continuation} of the first line. 
In other word, our method of computing sums of divergent series is analogous to the method of $\zeta$ regularization based on 
the Riemann's functional equality (the duality formula) \cite{Riemann}, which was derived also by analytical continuation. \\

For the concreteness let us give several explicit formulas for $\rho\geq 4$. For $\rho=4$ we have ($p=q^{-1}$)
\begin{align}
\begin{split}
S_4^{(1)}(q)&=\sum_{n=0}^\infty (-1)^n q^{n^2}
\\
&=1+\sum_{n=1}^\infty (-1)^n q^n \frac{(1+q)\cdots(1+q^{4n-3})}
{(1+q^3)\cdots(1+q^{4n-1})}
\\
&=1+\sum_{n=1}^\infty (-1)^n p^n \frac{(1+p)\cdots(1+p^{4n-3})}{(1+p^3)\cdots(1+p^{4n-1})}.
\end{split}
\end{align}
Interestingly the second and the third lines have the completely same form. 
Such phenomenon can be called {\it self-duality}, because $S_4^{(1)}(q)=S_4^{(1)}(q^{-1})$.  
In fact for $q=p=0.5$ both has the same value
\begin{align}
1-2^{-1}+2^{-4}-2^{-9}+\cdots&=1-2+2^4-2^9+\cdots
\nonumber \\
&=0.5605621040 \cdots.
\end{align}

Next for $\rho=5$ we have
\begin{align}
\begin{split}
S_5^{(1)}(q)&=\sum_{n=0}^\infty (-1)^n q^{n(3n-1)/2}
\\
&=1+\sum_{n=1}^\infty (-1)^n q^n \frac{(1+q)\cdots(1+q^{6n-5})}{(1+q^4)\cdots(1+q^{6n-2})}
\\
&=1+\sum_{n=1}^\infty (-1)^n p^{2n} \frac{(1+p)\cdots(1+p^{6n-5})}{(1+p^4)\cdots(1+p^{6n-2})},
\end{split}
\end{align}
which gives numerically for $q=p=0.5$
\begin{align}
\begin{split}
&1-2^{-1}+2^{-5}-2^{-12}+\cdots=0.5310060977\cdots,
\\
&1-2+2^5-2^{12}+\cdots=0.7181272344\cdots.
\end{split}
\end{align}

The first six results of $\rho=3, 4,\ \cdots, 8$ are shown in the following table.
\begin{center}
\begin{tabular}{cll}
$\rho$&$q=0.5$&$q=2\ (p=0.5)$\\
\hline
3&$ 0.610321518048\cdots$&$ 0.427525130255\cdots$\\
4&$ 0.560562104001\cdots$&$ 0.560562104001\cdots$\\
5&$ 0.531006097764\cdots$&$ 0.718127234453\cdots$\\
6&$ 0.515594486147\cdots$&$ 0.838420806786\cdots$\\
7&$ 0.507808685360\cdots$&$ 0.913124740197\cdots$\\
8&$ 0.503905773163\cdots$&$ 0.954897984724\cdots$
\end{tabular}
\end{center}
The limit behavior can be expected that as $\rho$ increases the former decreases to $0.5$ and the latter increases to $1.0$. 
In fact, we can show easily $S_\infty^{(1)}(q)=1-q$ for general $q<1$, and $S_\infty^{(1)}(q)=1$ for general $q>1$. 

Finally let us give results for another choice of $q=1/3$ and $q=3\ (p=1/3)$ for the trigonal number series, which are
\begin{align}
\begin{split}
&1-3^{-1}+3^{-3}-3^{-6}+\cdots=0.7023488270\cdots,
\\
&1-3+3^3-3^6+\cdots=0.4131192691\cdots,
\end{split}
\end{align}
which tells the cases $q=1/2$ and $q=2$ are not special. 


\subsection{Second group of problems -- Pochhammer symbol series}

The second Gauss problem is to find the value of
\begin{align}
\sum_{n=0}^\infty (q;q)_n=1+(1-q)+(1-q)(1-q^2)+\cdots
\end{align}
for $q=2$. However the sum of this series is {\it divergent} even for $q<1$, because $(q;q)_\infty=\prod_{n=1}^\infty (1-q^n)$ 
is finite. Therefore our answer to Gauss problem 2, is {\it divergent} against the expectation. 
More details and other examples of divergent series are discussed in a recent article by the present author \cite{SogoTriplicity}.

Let us consider next
\begin{align}
S_\kappa^{(2)}(q)&=\sum_{n=0}^\infty q^n (q; q^\kappa)_n
\\
&=1+q(1-q)+q^2(1-q)(1-q^{1+\kappa})+\cdots,
\nonumber
\end{align}
which can be written in a continued fraction form such that
\begin{align}
&S_\kappa^{(2)}(q)=\frac{d_0}{1+}\ \frac{d_1}{1+}\ \frac{d_2}{1+}\ \frac{d_3}{1+}\ \cdots,
\\
&d_0=1,\ d_1=q^2-q,\ d_2=q^{2+\kappa}-q^2,\ d_3=q^{2+2\kappa}-q^{1+\kappa},\ \cdots,\  
\nonumber \\
&d_{2m}=q^{2\kappa m+2-\kappa}-q^{\kappa m+2-\kappa},\ 
d_{2m+1}=q^{2\kappa m+2}-q^{\kappa m+1},
\label{Gauss4}
\end{align}
by using the Rogers formula \cite{Rogers}, which is explained next.

For the sake of completeness let us give the Rogers formula, which was found earlier by Muir \cite{Muir}. 
Suppose we have an equality
\begin{align}
c_0+c_1x+c_2x^2+c_3x^3+\cdots=
\frac{e_0}{1-}\ \frac{e_1x}{1-}\ \frac{e_2x}{1-}\ \frac{e_3x}{1-}\ \cdots,
\end{align}
then the coefficients $e$'s in the right hand side are given by
\begin{align}
\begin{split}
&e_0=\alpha_0,\ e_1=\frac{\alpha_1}{\alpha_0},\ e_2=\frac{\alpha_2}{\alpha_1\alpha_0},\ 
e_3=\frac{\alpha_3\alpha_0}{\alpha_2\alpha_1},\ \cdots,\\ 
&e_n=\frac{\alpha_n\alpha_{n-3}}{\alpha_{n-1}\alpha_{n-2}},\quad (n\geq 3),
\end{split}
\end{align}
where $\alpha$'s are determined in terms of $c$'s by
\begin{align}
\alpha_0=c_0,\ \alpha_1=c_1,\ 
\alpha_2=\left|\begin{array}{cc}
c_0&c_1\\
c_1&c_2\end{array}\right|,\ 
\alpha_3=\left|\begin{array}{cc}
c_1&c_2\\
c_2&c_3\end{array}\right|,\ \cdots
\end{align}
whose general expressions are given by ($m\geq 1$)
\begin{align}
\begin{split}
&\alpha_{2m}=\left|\begin{array}{cccc}
c_0&c_1&\cdots&c_{m}\\
c_1&c_2&\cdots&c_{m+1}\\
\cdots&\cdots&\cdots&\cdots\\
c_{m}&\cdots&\cdots&c_{2m}\end{array}\right|,\ 
\\
&\alpha_{2m+1}=\left|\begin{array}{cccc}
c_1&c_2&\cdots&c_{m+1}\\
c_2&c_3&\cdots&c_{m+2}\\
\cdots&\cdots&\cdots&\cdots\\
c_{m+1}&\cdots&\cdots&c_{2m+1}\end{array}\right|.
\end{split}
\end{align}
Such is the Rogers (or Muir-Rogers) formula. 

Applying these to the power series
\begin{align}
S_\kappa^{(2)}(x)=\sum_{n=0}^\infty (q; q^\kappa)_n x^n=\frac{e_0}{1-}\ \frac{e_1x}{1-}\ \frac{e_2x}{1-}\ \frac{e_3x}{1-}\ \cdots,
\end{align}
we obtain $(m\geq 1)$
\begin{align}
\begin{split}
&e_0=1,\ e_1=1-q,\ e_2=q(1-q^\kappa),\ e_3=q^{\kappa}(1-q^{\kappa+1}),\ \cdots,
\\
&e_{2m}=q^{\kappa m+1-\kappa}(1-q^{\kappa m}),\quad
e_{2m+1}=q^{\kappa m}(1-q^{\kappa m+1}),
\end{split}
\end{align}
which give \eqref{Gauss4} by setting $x=q$, which was announced before in \eqref{P_CF}. 

Now to proceed further, we have to borrow results of Ramanujan's continued fraction. Ramanujan wrote in his {\it Lost Notebook} 
\cite{RamanujanL, Andrews2} the following identity
\begin{align}
\frac{G(aq,\lambda q; b,q)}{G(a,\lambda; b,q)}=\frac{d_0}{1+}\ \frac{d_1}{1+}\ \frac{d_2}{1+}\ \frac{d_3}{1+}\ \cdots,
\end{align}
where he set
\begin{align}
\begin{split}
&G(a,\lambda;b,q)=1+\sum_{n=1}^\infty \frac{q^{n(n+1)/2}}{(q;q)_n}
\frac{(a+\lambda)\cdots(a+\lambda q^{n-1})}{(1+bq)\cdots(1+bq^n)},
\\
&G(aq,\lambda q;b,q)=1+\sum_{n=1}^\infty \frac{q^{n(n+3)/2}}{(q;q)_n}
\frac{(a+\lambda)\cdots(a+\lambda q^{n-1})}{(1+bq)\cdots(1+bq^n)},
\end{split}
\end{align}
then he derived
\begin{align}
\begin{split}
&d_0=1,\ d_1=aq+\lambda q,\ d_2=bq+\lambda q^2,\ d_3=aq^2+\lambda q^3,\  \cdots,
\\
&d_{2m}=bq^m+\lambda q^{2m},\quad 
d_{2m+1}=a q^{m+1}+\lambda q^{2m+1},\quad(m\geq 1),
\end{split}
\label{RamanujanCF}
\end{align}
an alternative proof of which is given in \cite{SogoTriplicity}.

Such is called the Ramanujan's continued fraction, which Ramanujan used reversely to obtain many identities by choosing 
parameters $a, \lambda, b$ and $q$ variously. For example, if we set $a=b=0,\ \lambda=1$ we obtain the Rogers-Ramanujan identity
\begin{align}
\frac{1+\sum_{n=1}^\infty \frac{q^{n(n+1)}}{(1-q)\cdots(1-q^n)}}{1+\sum_{n=1}^\infty \frac{q^{n^2}}{(1-q)\cdots(1-q^n)}}
=\frac{1}{1+}\ \frac{q}{1+}\ \frac{q^2}{1+}\ \frac{q^3}{1+}\ \cdots,
\end{align}
which was reported in 1919 by Ramanujan \cite{RamanujanC}, while Rogers found it earlier in 1894 \cite{RogersR} by using a different method.

If we compare \eqref{RamanujanCF} with \eqref{Gauss4}, we find the replacements in \eqref{RamanujanCF}
\begin{align}
q\rightarrow q^\kappa,\quad a\rightarrow -q^{1-\kappa},\quad \lambda\rightarrow q^{2-\kappa},\quad b\rightarrow -q^{2-\kappa},
\label{SecondGauss}
\end{align}
give the continued fraction \eqref{Gauss4}, because
\begin{align}
\begin{split}
&d_{2m}=bq^{m}+\lambda q^{2m}\rightarrow -q^{2-\kappa}q^{\kappa m}+q^{2-\kappa}q^{2\kappa m}
=q^{2\kappa m+2-\kappa}-q^{\kappa m+2-\kappa},
\\
&d_{2m+1}=aq^{m+1}+\lambda q^{2m+1}\rightarrow -q^{1-\kappa}q^{\kappa(m+1)}+q^{2-\kappa}q^{\kappa(2m+1)}
=q^{2\kappa m+2}-q^{\kappa m+1}.
\end{split}
\nonumber
\end{align}

For the later convenience, let us discuss here how Ramanujan's $G(a, \lambda; b, q)$ is obtained from Heine's $q$-hypergeometric series
\begin{align}
{}_2\phi_1(\alpha, \beta, \gamma; q, \tau)=\sum_{n=0}^\infty \frac{(\alpha; q)_n(\beta; q)_n}{(\gamma;q)_n(q;q)_n}\ \tau^n. 
\end{align}
If we take a limit of $\alpha\rightarrow\infty,\ \tau\rightarrow 0$ with a constraint $\alpha \tau=-qx$, we have
\begin{align}
\lim_{\tau\rightarrow 0}\ {}_2\phi_1(-qx/\tau, \beta, \gamma;q,\tau)=
\sum_{n=0}^\infty \frac{q^{n(n+1)/2}}{(q;q)_n}\frac{(\beta;q)_n}{(\gamma;q)_n}\ x^n,
\end{align}
which is equal to $G(a,\lambda; b, q)$ when $\beta=-\lambda/a,\ \gamma=-bq$ and $x=a$, that is,
\begin{align}
G(a,\lambda;b,q)&=\lim_{\tau\rightarrow 0}\ {}_2\phi_1(-aq/\tau, -\lambda/a, -bq; q, \tau).
\end{align}
Hereafter we wish to write
\begin{align}
\begin{split}
G_0&=G(a,\lambda;b,q)=\lim_{\tau\rightarrow 0}\ {}_2\phi_1(-aq/\tau, -\lambda/a, -bq; q, \tau)
\\
&=\lim_{\tau\rightarrow 0} \sum_{n=0}^\infty 
\frac{(-aq/\tau;q)_n (-\lambda/a;q)_n}{(-bq;q)_n(q;q)_n}\ \tau^n,
\end{split}
\\
\begin{split}
G_1&=G(aq,\lambda q; b, q)=\lim_{\tau\rightarrow 0}\ {}_2\phi_1(-aq^2/\tau, -\lambda/a, -bq; q, \tau)
\\
&=\lim_{\tau\rightarrow 0} \sum_{n=0}^\infty 
\frac{(-aq^2/\tau;q)_n (-\lambda/a;q)_n}{(-bq;q)_n(q;q)_n}\ \tau^n.
\end{split}
\end{align}

Next useful formulas are the transformation relations found by Heine \cite{Heine}, 
\begin{align}
{}_2\phi_1(\alpha,\beta,\gamma;q,\tau)
&=\frac{(\beta;q)_\infty (\alpha\tau;q)_\infty}{(\gamma;q)_\infty(\tau;q)_\infty}\ 
{}_2\phi_1(\gamma/\beta, \tau, \alpha\tau; q,\beta),
\label{Heine_1}
\\
&=\frac{(\gamma/\beta;q)_\infty(\beta \tau;q)_\infty}{(\gamma;q)_\infty(\tau;q)_\infty}\ 
{}_2\phi_1(\alpha\beta \tau/\gamma,\beta,\beta \tau; q, \gamma/\beta),
\label{Heine_2}
\\
&=\frac{(\alpha\beta\tau/\gamma;q)_\infty}{(\tau;q)_\infty}\ 
{}_2\phi_1(\gamma/\alpha, \gamma/\beta, \gamma; q, \alpha\beta\tau/\gamma),
\label{Heine_3}
\end{align}
which are (1.4.4), (1.4.5), (1.4.6) of the book by Gasper and Rahman \cite{GasperRahman}, 
and we set $(\alpha;q)_\infty=\prod_{n=0}^\infty (1-\alpha q^{n})$. 
Since ${}_2\phi_1$ has symmetry of $(\alpha, \beta)$ exchange, we have further
\begin{align}
{}_2\phi_1(\alpha,\beta,\gamma;q,\tau)
&=\frac{(\alpha;q)_\infty (\beta\tau;q)_\infty}{(\gamma;q)_\infty(\tau;q)_\infty} 
{}_2\phi_1(\gamma/\alpha, \tau, \beta\tau; q,\alpha),
\label{Heine_4}
\\
&=\frac{(\gamma/\alpha;q)_\infty(\alpha \tau;q)_\infty}{(\gamma;q)_\infty(\tau;q)_\infty}\ 
{}_2\phi_1(\alpha\beta \tau/\gamma,\alpha,\alpha \tau; q, \gamma/\alpha).
\label{Heine_5}
\end{align}
These five transformation relations can be applied to Ramanujan function $G_0$ and $G_1$. 

Let us take $\kappa=1$ case for the preliminary test, by using $a=-1,\ \lambda=q,\ b=-q$. Then we have
\begin{align}
\begin{split}
&G_0=G(-1,q;-q,q)=1+\sum_{n=1}^\infty (-1)^n\frac{q^{n(n+1)/2}}{(1-q^2)\cdots(1-q^{n+1})}
\\
&G_1=G(-q,q^2;-q,q)=1+\sum_{n=1}^\infty (-1)^n\frac{q^{n(n+3)/2}}{(1-q^2)\cdots(1-q^{n+1})},
\end{split}
\label{G0G1_1}
\end{align}
which are transformed by \eqref{Heine_2} into
\begin{align}
&G_0=\frac{(q;q)_\infty}{(q^2;q)_\infty}=1-q,\quad
G_1=\frac{(q;q)_\infty}{(q^2;q)_\infty}\sum_{n=0}^\infty (q;q)_n q^n=(1-q)\sum_{n=0}^\infty q^n(q;q)_n
\nonumber \\
&\Longrightarrow\quad 
S_1^{(2)}(q)=\frac{G_1}{G_0}=\sum_{n=0}^\infty q^n(q;q)_n
\end{align}
as it should be, where $(q;q)_\infty=(1-q)\cdot(q^2;q)_\infty$, and $(1; q)_0=1,\ (1; q)_n=0\ (n=1, 2, 3,\cdots)$ are used. 
The equality
\begin{align}
G_0=1+\sum_{n=1}^\infty (-1)^n\frac{q^{n(n+1)/2}}{(1-q^2)\cdots(1-q^{n+1})}=1-q,
\end{align}
is a strange but remarkable result which is also checked numerically for $0\leq q<1$. 

In concluding the case $\kappa=1$, we have
\begin{align}
\begin{split}
S_1^{(2)}(q)&=\sum_{n=0}^\infty q^n(q;q)_n
=\frac{1}{1+}\ \frac{q^2-q}{1+}\ \frac{q^3-q^2}{1+}\ \frac{q^4-q^2}{1+}\ \frac{q^5-q^3}{1+}\ \cdots
\\
&=\frac{G_1}{1-q}=\sum_{n=0}^\infty (-1)^n \frac{q^{n(n+3)/2}}{(q;q)_{n+1}},
\end{split}
\label{S12} 
\end{align}
where \eqref{G0G1_1} is used for $G_1$, and $(1-q)(1-q^2)\cdots(1-q^{n+1})=(q;q)_{n+1}$. 
We can rewrite \eqref{S12} by substituting $q=p^{-1}$, such that
\begin{align}
\begin{split}
S_1^{(2)}(q)&=\sum_{n=0}^\infty q^n(q;q)_n
\\
&=\sum_{n=0}^\infty (-1)^n \frac{q^{n(n+3)/2}}{(q;q)_{n+1}}
\\
&=-\left(\sum_{n=0}^\infty \frac{p}{(1-p)(1-p^2)\cdots(1-p^{n+1})}\right),
\end{split}
\end{align}
which can be used to define $S_1^{(2)}(q)$ for $q>1$ by the right hand side series for $p<1$. 
However since it is still divergent, we have to conclude the series $S_1^{(2)}(q)$ for $q>1$ is {\it divergent}, 
although it is convergent for $q<1$. 

For the next $\kappa=2$ case, the parameters are $q \rightarrow q^2$ and $a=-q^{-1},\ \lambda=1,\ b=-1$, and we obtain
\begin{align}
\begin{split}
S_2^{(2)}(q)&=\sum_{n=0}^\infty q^n(q;q^2)_n
=\frac{1}{1+}\ \frac{q^2-q}{1+}\ \frac{q^4-q^2}{1+}\ \frac{q^6-q^3}{1+}\ \frac{q^8-q^4}{1+}\ \cdots
\\
&=\sum_{n=0}^\infty (-1)^n\frac{q^{n(n+1)}}{(q;q^2)_{n+1}},
\end{split}
\label{GaussP3}
\end{align}
which was given also by Andrews (see p.181 of \cite{Andrews2}). The second line is derived as follows. Firstly
\begin{align}
G_0&=G(-q^{-1}, 1; -1, q^2)=\lim_{\tau\rightarrow 0}\sum_{n=0}^\infty 
\frac{(q/\tau;q^2)_n(q;q^2)_n}{(q^2;q^2)_n(q^2;q^2)_n}\ \tau^n
\nonumber \\
&=\lim_{\tau\rightarrow 0}\frac{(q;q^2)_\infty(q\tau;q^2)_\infty}{(q^2;q^2)_\infty(\tau;q^2)_\infty}
\sum_{n=0}^\infty \frac{(q;q^2)_n(1;q^2)_n}{(q^2;q^2)_n(q\tau;q^2)_n}\ q^n
\nonumber \\
&=\frac{(q;q^2)_\infty}{(q^2;q^2)_\infty},
\end{align}
where \eqref{Heine_2} and $(1;q^2)_n=0\ (n\geq 1)$ are used. Nextly 
\begin{align}
G_1&=G(-q, q^2; -1,q^2)=
\lim_{\tau\rightarrow 0}\sum_{n=0}^\infty \frac{(q^3/\tau;q^2)_n(q;q^2)_n}{(q^2;q^2)_n(q^2;q^2)_n}\ \tau^n
\nonumber \\
&=\lim_{\tau\rightarrow 0}\frac{(q^3;q^2)_\infty(\tau/q;q^2)_\infty}{(q^2;q^2)_\infty(\tau;q^2)_\infty}
\sum_{n=0}^\infty \frac{(q^3/\tau;q^2)_n(q^2;q^2)_n}{(q^3;q^2)_n(q^2;q^2)_n}(\tau/q)^n
\nonumber \\
&=\frac{(q^3;q^2)_\infty}{(q^2;q^2)_\infty}\sum_{n=0}^\infty (-1)^n\frac{q^{n(n+1)}}{(q^3;q^2)_n},
\end{align}
where \eqref{Heine_5} is used. Therefore we obtain
\begin{align}
S_2^{(2)}(q)=\frac{G_1}{G_0}=\sum_{n=0}^\infty (-1)^n\frac{q^{n(n+1)}}{(q;q^2)_{n+1}}
\end{align}
where $(q;q^2)_\infty=(1-q)\cdot(q^3;q^2)_\infty$, and $(q;q^2)_{n+1}=(1-q)\cdot(q^3;q^2)_n$ are used.

We can rewrite \eqref{GaussP3} by setting $q=p^{-1}$
\begin{align}
S_2^{(2)}(q)&=\sum_{n=0}^\infty q^n(q;q^2)_n
\nonumber \\
&=\sum_{n=0}^\infty (-1)^n\frac{q^{n(n+1)}}{(q; q^2)_{n+1}}
=-\left(\sum_{n=0}^\infty \frac{p^{n+1}}{(p; p^2)_{n+1}}\right),
\end{align}
which enables us to define the value of $S_2^{(2)}(q)$ for $q>1$ by the right hand side series which is convergent for $p<1$. 
For example when $p=0.5\ (q=2)$ we obtain
\begin{align}
\begin{split}
S_2^{(2)}(2)&=\sum_{n=0}^\infty 2^n(2;2^2)_n=1-2+4\cdot7-8\cdot 7\cdot 31+\cdots
\\
&=-\left(\sum_{n=0}^\infty \frac{2^{-(n+1)}}{(2^{-1};2^{-2})_{n+1}}\right)
=-2.1639450388\cdots.
\end{split}
\label{GaussP3R}
\end{align}
The author reported this result recently as the answer to {\it Gauss problem 3} \cite{SogoTriplicity}.\\

Let us consider arbitrary $\kappa$ case, which corresponds to $q\rightarrow q^\kappa$ and $a=-q^{1-\kappa},\ \lambda=q^{2-\kappa},\ 
b=-q^{2-\kappa}$ by using \eqref{SecondGauss}, which give $\alpha=-aq^\kappa/\tau=q/\tau,\ \beta=-\lambda/a=q,\ \gamma=-bq^\kappa=q^2$. 
Therefore we have
\begin{align}
G_0&=G(a, \lambda; b, q^\kappa)=\sum_{n=0}^\infty (-1)^n 
\frac{q^{n[\kappa(n-1)+2]/2}}{(q^\kappa;q^\kappa)_n}\ \frac{(q;q^\kappa)_n}{(q^2;q^\kappa)_n},
\nonumber \\
&=\lim_{\tau\rightarrow 0}
\sum_{n=0}^\infty \frac{(q/\tau;q^\kappa)_n(q;q^\kappa)_n}{(q^\kappa;q^\kappa)_n(q^2;q^\kappa)_n}\ \tau^n
=\frac{(q;q^\kappa)_\infty}{(q^2;q^\kappa)_\infty},
\\
G_1&=G(aq^\kappa,\lambda q^\kappa;b, q^\kappa)
=\sum_{n=0}^\infty (-1)^n \frac{q^{n[\kappa(n-1)+4]/2}}{(q^\kappa;q^\kappa)_n}\ \frac{(q;q^\kappa)_n}{(q^2;q^\kappa)_n}
\nonumber \\
&=\lim_{\tau\rightarrow 0}\sum_{n=0}^\infty 
\frac{(q^{1+\kappa}/\tau, q^\kappa)_n(q;q^\kappa)_n}{(q^\kappa;q^\kappa)_n(q^2;q^\kappa)_n}\ \tau^n
\nonumber \\
&=\frac{(q^{1+\kappa};q^\kappa)_\infty}{(q^2;q^\kappa)_\infty}\sum_{n=0}^\infty (-1)^n 
\frac{q^{n[\kappa(n-1)+4]/2}}{(q^{1+\kappa};q^\kappa)_n},
\end{align}
where \eqref{Heine_2} for $G_0$, \eqref{Heine_5} for $G_1$ is applied respectively.
These results coincide with those for $\kappa=1$ and $\kappa=2$ obtained already. 

For concreteness let us give results of $\kappa=3$ case explicitly. They are given by
\begin{align}
\begin{split}
G_0&=G(a,\lambda;b,q^3)=\frac{(q;q^3)_\infty}{(q^2;q^3)_\infty},
\\
G_1&=G(aq^3,\lambda q^3; b, q^3)=
\frac{(q^4;q^3)_\infty}{(q^2;q^3)_\infty}\sum_{n=0}^\infty (-1)^n\frac{q^{n(3n+1)/2}}{(q^4;q^3)_n},
\end{split}
\label{T0}
\end{align}
therefore we have
\begin{align}
\begin{split}
S_3^{(2)}(q)&=\sum_{n=0}^\infty q^n(q;q^3)_n
\\
&=\frac{G_1}{G_0}=\sum_{n=0}^\infty (-1)^n\frac{q^{n(3n+1)/2}}{(q;q^3)_{n+1}}
\\
&=-\left(\sum_{n=0}^\infty \frac{p^{2n+1}}{(p;p^3)_{n+1}}\right),\quad (p=q^{-1})
\end{split}
\end{align}
where $(q;q^3)_\infty=(1-q)\cdot(q^4;q^3)_\infty$ and $(1-q)\cdot(q^4;q^3)_n=(q;q^3)_{n+1}$ are used.

Again such relation enables us to define the value of $S_3^{(2)}(q)$ for $q>1$. 
For example when $p=0.5\ (q=2)$ we obtain
\begin{align}
\begin{split}
S_3^{(2)}(2)&=\sum_{n=0}^\infty 2^n(2;2^3)_n=1-2+4\cdot15-8\cdot 15\cdot 127+\cdots
\\
&=-\left(\sum_{n=0}^\infty \frac{2^{-(2n+1)}}{(2^{-1};2^{-3})_{n+1}}\right)
=-1.3562780680\cdots.
\end{split}
\label{GaussP3R3}
\end{align}

In conclusion the general formula for arbitrary $\kappa$ can be summarized by the following theorem.\\

\noindent
{\bf Theorem 2}\quad For $\kappa=1,2,3,\cdots$, we have
\begin{align}
\begin{split}
S_\kappa^{(2)}(q)&=\sum_{n=0}^\infty q^n(q;q^\kappa)_n
\\
&=\sum_{n=0}^\infty (-1)^n \frac{q^{n(\kappa n+4-\kappa)/2}}{(q; q^\kappa)_{n+1}}
\\
&=-\left(\sum_{n=0}^\infty \frac{p^{(\kappa-1) n+1}}{(p; p^\kappa)_{n+1}}\right), \quad(p=q^{-1}).
\end{split}
\label{GF_2}
\end{align}

Let us consider again the parameter $q$ is a complex number, and $\kappa\geq 2$. 
The first line and the second line are convergent for $|q|<1$, and equivalent with each other because both have 
the same continued fraction expression. 
Now the second line and the third line are the same formula because the latter is just a rewriting of the former 
by substitution $p=q^{-1}$. 
Since the third line is convergent for $|p|<1$ {\it i.e.}, $|q|>1$, 
we can say that it is {\it analytical continuation} of the first line. 
Therefore the method of computing sums of divergent series in {\bf theorem 2} is the same as in {\bf theorem 1}, that is, 
again based on the concept of analytical continuation. \\

The first six results of $\kappa=1, 2,\ \cdots, 6$ for $q=1/2$ and $2$ are shown in the following table. 

\begin{center}
\begin{tabular}{cll}
$\kappa$& $q=0.5$ & $q=2\ (p=0.5)$\\
\hline
1&$1.422423809826\cdots$& $-\infty$\\
2&$1.464859485508\cdots$&$-2.163945038886\cdots$\\
3&$1.483398918062\cdots$&$-1.356278068046\cdots$\\
4&$1.491943374305\cdots$&$-1.147501794618\cdots$\\
5&$1.496032715309\cdots$&$-1.067726939580\cdots$\\
6&$1.498031616225\cdots$&$-1.032512189123\cdots$
\end{tabular}
\end{center}
which suggests that $S_\kappa^{(2)}(1/2)\rightarrow 1.5$ and $S_\kappa^{(2)}(2)\rightarrow -1.0$ when $\kappa\rightarrow\infty$. 
In fact, we can show that $S_\infty^{(2)}(q)=1+q$ for $q<1$, and $S_\infty^{(2)}(q)=-p/(1-p)=1/(1-q)$ for $q>1$. 


\section{Concluding remarks}
\setcounter{equation}{0}

Motivated by Entry 7 of Gauss's diary, two groups of series $S_\rho^{(1)}(q)$ and $S_\kappa^{(2)}(q)$,
\begin{align}
&S_\rho^{(1)}(q)=\sum_{n=0}^\infty q^{n[(\rho-2)n-(\rho-4)]/2},\quad(\rho=3,4,5,\cdots)
\\
&S_\kappa^{(2)}(q)=\sum_{n=0}^\infty q^n(q;q^\kappa)_n,\quad(\kappa=1,2,3,\cdots)
\end{align}
are considered. They are convergent for $q<1$, but are divergent for $q>1$. 
Our idea is to find another equivalent expression for the series and replace the original series by them. 
Such series are found as follows, 
\begin{align}
&S_\rho^{(1)}(q)
=\sum_{n=0}^\infty (-1)^n q^n\ \frac{(-q; q^{2\rho-4})_n}{(-q^{\rho-1}; q^{2\rho-4})_n}
=\sum_{n=0}^\infty (-1)^n p^{(\rho-3)n} \frac{(-p; p^{2\rho-4})_n}{(-p^{\rho-1}; p^{2\rho-4})_n},
\\
&S_\kappa^{(2)}(q)
=\sum_{n=0}^\infty (-1)^n \frac{q^{n(\kappa n+4-\kappa)/2}}{(q; q^\kappa)_{n+1}}
=-\left(\sum_{n=0}^\infty \frac{p^{(\kappa-1) n+1}}{(p; p^\kappa)_{n+1}}\right),
\end{align}
which are \eqref{GF_1} of {\bf Theorem 1}, and \eqref{GF_2} of {\bf Theorem 2}. 
In the right hand side, we set $p=q^{-1}<1$ when $q>1$. 
This idea of introducing $p=q^{-1}$ was inspired by Slater's book \cite{Slater}, specifically her discussions in \S 3.2 
on convergence property of basic hypergeometric series. 

Our method of computing sums of divergent series is based on the concept of {\it analytical continuation}, 
which is described in the previous section after the statements of {\bf theorem 1} and {\bf theorem 2}. 
Therefore the uniqueness of the results can be attributed to the uniqueness theorem of analytical continuation or 
{\it identity theorem}. 

Our logic of computing sums of divergent series can be interpreted in another way. 
Euler used \cite{HardyD, Bourbaki} 
\begin{align}
\begin{split}
&1-x+x^2-x^3+\cdots=\frac{1}{1+x}=x^{-1}-x^{-2}+x^{-3}-\cdots
\\
&\Longleftrightarrow\quad 
\sum_{n=-\infty}^\infty (-1)^n x^n=0,
\\
&1+x+x^2+x^3+\cdots=\frac{1}{1-x}=-x^{-1}-x^{-2}-x^{-3}-\cdots
\\
&\Longleftrightarrow\quad 
\sum_{n=-\infty}^\infty x^n=0,
\end{split}
\label{Euler}
\end{align}
when he wrote $1-2+4-8+\cdots=1/3$, and $1+2+4+8+\cdots=-1$. 
Since such trick is almost the same with the present method, our method may be called {\it modified Euler's method}, 
because the left hand side series for $|x|<1$ and the right hand side series for $|x|>1$ in \eqref{Euler} are 
analytical continuations of the central function $(1+x)^{-1}$ or $(1-x)^{-1}$. 
And the curious equalities in \eqref{Euler} should be rewritten, in modern terminology, by
\begin{align}
\begin{split}
&\sum_{n=-\infty}^\infty (-1)^n\ e^{2\pi i ny}=\sum_{m=-\infty}^\infty \delta\left(y-(m+\frac{1}{2})\right),
\\
&\sum_{n=-\infty}^\infty e^{2\pi iny}=\sum_{m=-\infty}^\infty \delta(y-m),
\end{split}
\label{Last}
\end{align}
where we set $x=e^{2\pi iy}$ and $\delta(y)$ is Dirac's delta function. The right hand sides of \eqref{Last} vanish 
except $y$ is half-integer or integer respectively. \\




%
\end{document}